\documentclass[11pt]{amsart}
\usepackage{amsmath,amssymb,amsfonts,url,mathptmx,tikz,geometry}
\usepackage{color}
\usepackage{appendix}
\numberwithin{equation}{section}

\newtheorem{thm}{Theorem}[section]
\newtheorem{lem}{Lemma}[section]
\newtheorem{rem}{Remark}[section]

\newtheorem{theoremA}{Theorem}

\usepackage{amsthm,a0size,bm,indentfirst,wasysym}

\usepackage{amsmath}

\allowdisplaybreaks[4]

\begin{document}

\title[Uniqueness of blow-up solution]{A simple proof of the  Uniqueness of blow-up solutions of mean field equations}

\author[L. Wu]{Lina Wu}
\address{Lina Wu, School of Mathematics and Statistics \\
	Beijing Jiaotong University \\
	Beijing, 100044, China}
\email{lnwu@bjtu.edu.cn}

\author[W. zou]{Wenming Zou}
\address{Wenming Zou, Department of Mathematical Sciences \\ 
        Tsinghua University \\ 
        Beijing, 100084, China}
\email{zou-wm@mail.tsinghua.edu.cn}

\date{\today}
	

\begin{abstract}
For a regular mean field equation defined on a compact Riemann surface, an important work of Bartolucci-Jevnikar-Lee-Yang \cite{bart-4} proved a uniqueness theorem for blow-up solutions under non-degeneracy assumptions. However, the proof is highly nontrivial and challenging to read. In this article, we not only provide a simple proof for the regular equation but also extend our proof to the case of singular equations with negative singular poles. Our proof supplements what is not written in a recent outstanding work by Bartolucci-Yang-Zhang \cite{byz-1}. 
\end{abstract}

	
\maketitle
	
\section{Introduction}

It is widely recognized that the following mean field equation with singular sources:
\begin{equation}\label{m-equ}
\Delta_g v+\rho\bigg(\frac{he^v}{\int_M h e^v{\rm d}\mu}-\frac{1}{vol_g(M)}\bigg)=\sum_{j=1}^N 4\pi \alpha_j \left(\delta_{q_j}-\frac{1}{vol_g(M)}\right) \quad {\rm in} \ \; M,
\end{equation}
plays a significant role in various fields, including conformal geometry, Electroweak and Self-Dual Chern-Simons vortices (\cite{ambjorn, spruck-yang,taran-1,taran-2,y-yang}), the statistical mechanics description of turbulent Euler flows, plasmas, and self-gravitating systems (\cite{bart-5,caglioti-2,SO,TuY,wolan}), the theory of hyperelliptic
curves and modular forms (\cite{chai}) and CMC immersions in hyperbolic 3-manifolds (\cite{taran-4}).
Here $(M,g)$ is a Riemann surface with metric $g$, $\Delta_g$ is the Laplace-Beltrami operator ($-\Delta_g\ge 0$), $h$ is a positive function on $M$ satisfying $h\in C^5(M)$ , $q_1,\cdots,q_N$ are distinct points on $M$, $\rho>0,\alpha_j>-1$ are constants, $\delta_{q_j}$ is the Dirac measure at $q_j\in M$, $vol_g(M)$ is the volume of $(M,g)$.  For simplicity, we assume $vol_g(M)=1$ throughout this article. In a recent outstanding work of Bartolucci-Yang-Zhang, the authors presented an all around type result on the uniqueness of blow-up solutions of the main equation. 

\medskip

To reformulate the main equation in an equivalent form, we employ the standard Green's function $G_M(z,p)$, which is uniquely determined (\cite{Aub}) as the solution to the following boundary value problem:
\begin{equation*}
\begin{cases}
-\Delta_g G_M(x,p)=\delta_p-1\quad {\rm in}\ \; M
\\
\\
\int_{M}G_M(x,p){\rm d}\mu=0.
\end{cases}
\end{equation*}
In local isothermal coordinates $z=z(x)$ centered at $p$ (with $z(p)=0$), the function
$G(z,0)=G_M(x(z),p)$ admits the following decomposition:
$$G(z,0)=-\frac 1{2\pi }\log |z|+R(z,0),$$
where $R(z,0)$ represents the regular part in local coordinates. Using $G_M(x,p)$, we rewrite (\ref{m-equ}) in the following equivalent form:
\begin{equation}\label{r-equ}
    \Delta_g  w+\rho\bigg(\frac{He^w}{\int_M H e^w{\rm d}\mu}-1\bigg)=0 \quad {\rm in}\ \; M,
\end{equation}
where the functions $w$ and $H$ are defined as:
\begin{equation}\label{r-sol}
	w(x)=v(x)+4\pi \sum_{j=1}^N \alpha_j G_M(x,q_j),
\end{equation}
and
\begin{equation*}
H(x)=h(x)\prod_{j=1}^N e^{-4\pi\alpha_j G_M(x,q_j)}.
\end{equation*}
In local coordinate $z$ centered at each $q_j$ (with $z(q_j)=0$), the function $H$ exhibits the following behavior for $|z|\ll 1$
\begin{equation*}
H(z)=h_j(z)|z|^{2\alpha_j},\quad  1\leq j\leq N,
\end{equation*}
with each $h_j(z)$ is a strictly positive smooth function in a neighborhood of $z=0$.

A sequence $\{v_k\}$ is called a sequence of \textit{bubbling solutions} of (\ref{m-equ}) if the corresponding function $w_k$ defined by (\ref{r-sol}) satisfy $\parallel w_k\parallel \to +\infty$ as $k\to+\infty$. In such cases, it is established in \cite{BM3,BT,li-cmp} that $w_k$ exhibits blow-up behavior at a finite number of points. Let $\{p_1,\cdots,p_m\}$ denote the set of blow-up points. This implies there exist $m$ sequences of points $p_{k,1},\cdots,p_{k,m}$ such that,after possibly along a sequences, 
\[w_k(p_{k,j})\to +\infty\quad\text{and }  \quad  p_{k,j}\to p_j \, (j=1,\cdots,m)  \quad  \text{as}  \quad  k\to +\infty.\]
Each blow-up point $p$ carries a singular source with strength $4\pi \alpha_p$, where $\alpha_p=0$ if $p$ is regular. We denote the strengths at $p_1,\cdots,p_m$ by $\alpha_1,\cdots,\alpha_m$ respectively, and we let the first $\tau$ of them be singular blow-up points:
\begin{equation}\label{largest-s}
\alpha_i>-1,\quad \alpha_i\not \in \mathbb N, \quad 1\le i\le \tau,\quad \alpha_{\tau+1}=\cdots=\alpha_m=0.
\end{equation}
Here, $\mathbb N\cup \{0\}$ denotes the set of natural numbers including $0$ and
$\alpha_M:=\max\{\alpha_1,...,\alpha_m\}$.
We say that the singular source located at $p$ is \emph{non-quantized} if $\alpha_p$ is not a positive integer. Here and in the rest of this work, we assume that all the singular sources, as far as they happen to be blow-up points, are non-quantized.  

Equation (\ref{r-equ}) possesses a natural variational structure, serving as the Euler-Lagrange equation for the energy functional:
$$I_{\rho}(w)=\frac 12 \int_M |\nabla w|^2+\rho\int_M w-\rho \log \int_M He^w,\quad w\in H^1(M).$$
Since adding a constant to any solution of (\ref{r-equ}) still yields a solution, there is no loss of generality in restricting $I_{\rho}$ on the subspace of $H^1(M)$ functions with zero mean. A complete discussion of the variational structure of (\ref{r-equ}) can be found in \cite{mal-ruiz}.

\smallskip

To state the main results, we introduce the following notation:
\begin{equation*}
G_j^*(x)=8\pi (1+\alpha_j)R_M(x,p_j)+8\pi \sum_{l,l\neq j}(1+\alpha_l)G_M(x,p_l),
\end{equation*}
where $R_M(\cdot, p_j)$ denotes the regular part of $G_M(\cdot, p_j)$.
\begin{equation}\label{L-p}
L(\mathbf{p})=\sum_{j\in I_1}[\Delta \log h(p_j)+\rho_*-N^*-2K(p_j)]h(p_j)^{\frac{1}{1+\alpha_M}}e^{\frac{G_j^*(p_j)}{1+\alpha_M}},
\end{equation}
\begin{equation}
\label{a-note}
\begin{cases} 
\alpha_M=\max_{i}\alpha_i,\quad
I_1=\{i\in \{1,\cdots,m\};\quad \alpha_i=\alpha_M\}, \\
\\
\rho_*=8\pi\sum_{j=1}^m(1+\alpha_j),\quad  N^*=4\pi\sum_{j=1}^m\alpha_j. 
\end{cases}
\end{equation}
In other words, $L(\mathbf{p})$ evaluates the weighted sum over all indices $j$ where $\alpha_j$ attains its maximal value $ \alpha_M$ within the set $\{\alpha_1,\cdots,\alpha_m\}$. As mentioned above, whenever some blow-up points are non-quantized singular sources while some others are regular, we maintain the following convention without loss of generality:
\begin{itemize}
    \item $p_1,\cdots,p_{\tau}$ are singular sources (with $\alpha_i\notin \mathbb{N}$).
    \item $p_{\tau+1},\cdots,p_m$ are regular points (with $\alpha_i=0$).
\end{itemize}
In this case
for $(x_{\tau+1},\cdots,x_m)\in M\times\cdots \times M$ we define
\begin{equation*}
\begin{aligned}
  f^*(x_{\tau+1},\cdots,x_m)
  =&\sum_{j=\tau+1}^m\big[\log h(x_j)+4\pi R_M(x_j,x_j)\big]\\
&+4\pi \sum_{l\neq j}^{\tau +1,\cdots,m}G_M(x_l,x_j)
+\sum_{j=\tau+1}^m\sum_{i=1}^{\tau}8\pi(1+\alpha_i)G_M(x_j,p_i). 
\end{aligned}
\end{equation*}

It is well known (\cite{chen-lin-sharp}) that  $(p_{\tau+1},\cdots,p_m)$ is a critical point of $f^*$. In a recent article of Bartolucci-et-al \cite{byz-1}, the authors proved the uniqueness of bubbling solutions for the case $\alpha_M>0$. They also indicated that by using their approach together with that in \cite{bart-4}, the case $\alpha_M\le 0$ also follows. Since the proof in \cite{bart-4} (as well as the one in \cite{byz-1}) is extremely technical, in this article, we provide a short proof for the case $\alpha_M\le 0$ using the approach in \cite{byz-1}.

Obviously, if $\alpha_M>0$, the set of blow-up points contains at least one positive singular source. When $\alpha_M\le 0$, the blow-up set either consists of a mixture of regular blow-up points and negative singular sources ($\alpha_M=0$), or only contains negative singular sources $(\alpha_M<0)$. The result for the case $\alpha_M>0$ is already established in \cite{byz-1}.

\begin{theoremA}

    ($\alpha_M>0$)  Let $v_k^{(1)}$ and $v_k^{(2)}$ be two sequences of bubbling solutions of {\rm (\ref{m-equ})}  with the same $\rho_k$: $\rho_k^{(1)}=\rho_k=\rho_k^{(2)}$ and the same blow-up set: $\{p_1,\cdots,p_m\}$. Suppose $(\alpha_1,\cdots,\alpha_N)$ satisfies {\rm (\ref{largest-s})}, $\alpha_M>0$, $L(\mathbf{p})\neq 0$ and, as far as $\tau<m$,  $\det \big(D^2f^*(p_{\tau+1},\cdots,p_m)\big)\neq 0$. Then there exists $n_0>1$ such that
 $v_k^{(1)}=v_k^{(2)}$ for all $k\ge n_0$.
 
\end{theoremA}

In this context, $D^2f^*$ represents the Hessian tensor field on $M$, and when evaluating $\det(D^2f^*)$ at $p_j$ for $j=\tau+1,\cdots,m$, we use $\phi_j$ to denote the conformal factor, which satisfies 
\begin{equation*}
\Delta_g=e^{-\phi_j}\Delta,\quad \phi_j(0)=|\nabla \phi_j(0)|=0, \quad e^{-\phi_j}\Delta \phi_j=-2K,
\end{equation*}
for $j=\tau+1,\cdots,m$. Additionally, $h_j$ is interpreted as $h_je^{\phi_j}$ in a small neighborhood around $p_j$ (for $j=\tau+1,\cdots,m$).
Therefore, if $\tau=m$, meaning all blow-up points are singular, the only relevant assumption is that $L(\mathbf{p})\neq 0$. 

To describe the result for $\alpha_M\le 0$ (where the set of blow-up points only consists of regular points and negative sources), we shall introduce new quantities. We use $B(q,r)$ to denote the geodesic ball of radius $r$ centered at $q\in M$, and $\Omega(q,r)$ to represent the pre-image of the Euclidean ball of radius $r$, $B_r(q)\subset \mathbb R^2$, in a suitably defined isothermal coordinate system. We fix a family of open, mutually disjoint sets $M_j$, whose union's closure is the whole $M$. Now, we define
$$D(\mathbf{p})=\lim_{r\to 0}\sum_{j=1}^mh(q_j)e^{G_j^*(q_j)}\bigg (
\int_{M_j\setminus \Omega(q_j,r_j)}\rho^*
e^{\Phi_j(x,\mathbf{q})}d\mu(x)-\frac{\pi}{1+\alpha_j}r_j^{-2-2\alpha_j}\bigg ),$$
where $M_1=M$ if $m=1$, $\displaystyle{r_j=r \left(\frac{h(q_j)e^{G_j^*(q_j)}}{8(1+\alpha_j)^2}\right)^{\frac{1}{2\alpha_j+2}}}$ and
\begin{align*}\Phi_j(x,\mathbf{q})=&\sum_{l=1}^m8\pi(1+\alpha_l)G(x,q_l)-G_j^*(q_j)\\
&+\log h_j(x)-\log h_j(q_j)+4\pi\alpha_j(R_M (x,q_j)-G_M(x,q_j)),
\end{align*}
where we remark that in local isothermal coordinates centered at $q_j$ (thus $0=z(q_j)$), we have 
\[4\pi\alpha_j(R(z,0)-G(z,0))=2\alpha_j\log|z|.\] 
It is crucial to emphasize that in this case, as long as $\alpha_j\neq 0$, $\alpha_j\in (-1,0)$, and specifically, the non-integrable terms in the definition of $D(\mathbf{p})$ cancel out, ensuring that the limit exists. Therefore, we obtain the following:

\begin{thm}\label{mainly-case-2}
Under the same assumptions of Theorem A, but with the condition that now $\alpha_M\le 0$, in any one of the following situations:
\begin{enumerate}
    \item $\alpha_M=0$,\, $L(\mathbf{p})\neq 0$,\, $\det \big(D^2f^*(p_{\tau+1},\cdots,p_m)\big)\neq 0$
    \item $\alpha_M=0$, $L(\mathbf{p})=0$,  $D(\mathbf{p})\neq 0$, $\det \big(D^2f^*(p_{\tau+1},\cdots,p_m)\big)\neq 0$
    \item $\alpha_M<0$, $D(\mathbf{p})\neq 0$.
\end{enumerate} the same uniqueness property as in Theorem A holds.
\end{thm}

Theorem \ref{mainly-case-2} was initially stated in \cite{byz-1} without proof, where the authors left it to the readers to complete the proof themselves. However, they referred to \cite{bart-4} for the case where $\alpha_M=0$ and to \cite{wu-zhang-ccm} for the case $0<\alpha_M<1$.  The proofs in \cite{bart-4,wu-zhang-ccm} are quite complex. In this article, we utilize the ideas from \cite{byz-1} to provide a simpler proof of Theorem \ref{mainly-case-2}. The first two cases in Theorem \ref{mainly-case-2} cover the situations where all blow-up points are regular and where blow-up points include both regular points and negative singular poles. Specifically, we not only provide a simple proof for the regular equation but also extend our proof to the case of singular equations with negative singular poles.

Similar conclusions can also be drawn for the corresponding Dirichlet problem. The same proof also applies to Theorem 1.4 in \cite{byz-1}.

Although we primarily use the approach from \cite{byz-1}, there are some parts where their refined estimates are not actually sufficiently carried out. We need to improve the estimates for the bubbling solutions around the regular blow-up points and to obtain a uniform precision in estimates around all possible blow-up points. Our aim in this work is to supplement the proof of Theorem \ref{mainly-case-2}.

It should be noted that the proof of the uniqueness theorem mentioned above involves three parameters, namely $b_0$, $b_1$, and $b_2$, as introduced in Section 4. We aim to show that all three parameters are equal to zero. The proof of $b_0=0$ is more intricate when $\alpha_M\le 0$, as the dominant term of $b_0$ tends to zero at a relatively faster rate, which requires more refined estimates compared to the case when $\alpha_M> 0$. This aspect constitutes a distinctive feature of our paper. On the other hand, the proofs for $b_1=0$ and $b_2=0$ become more complex when $\alpha_M>0$, more precisely, when $\alpha_M>1$. At this stage, it is necessary to establish the precise vanishing estimates of $b_0^k$, as detailed in \cite{byz-1}. 

\medskip

To introduce the subsequent content more concisely and accurately, we present the following notational conventions:
\begin{itemize}
    \item $B_\tau=B_\tau(0)$ refers to a ball centered at the origin in local isothermal coordinates $y\in B_\tau$. 
    \item If $B_\tau$ is any such ball centered at some point $0=y(p), p\in M$,  its pre-image on $M$ is denoted $\Omega(p,\tau)\subset M$.
    \item $B(p,\tau)\subset M$ always denote a geodesic ball.
    \item In various estimates, a small positive number $\epsilon_0$ is used, and its value may vary across different contexts.
\end{itemize}


The organization of this paper is as follows. In section 2, we gather the preliminary estimates required for the proof of the main theorem. In section 3, we derive more refined estimates for the local equations. Finally, in section 4, we present the proofs of the main theorem.



\section{Preliminary Estimates}\label{preliminary}

Given that the proof of the main theorems involves intricate analysis, in this section, we present several estimates established in \cite{BCLT,BM3,BT,BT-2,chen-lin-sharp,chen-lin,li-cmp,zhang1,zhang2}.

Consider a sequence of solutions $w_k$ of \eqref{r-equ} with $\rho =\rho_k$, where $w_k$ exhibits blow-up behavior at $m$ points $\{p_1, \cdots,p_m\}$. To characterize the bubbling profile of $w_k$ near each $p_j$, we define the normalized function:
\begin{equation*}
	u_k=w_k-\log\bigg(\int_M He^{w_k}{\rm d}\mu \bigg),
\end{equation*}
which satisfies the normalization condition
$$\int_{M}He^{u_k}{\rm d}\mu=1,$$ 
The corresponding equation for  $u_k$ becomes
\begin{equation*}
\Delta_g u_k+\rho_k(He^{u_k}-1)=0\quad {\rm in} \ \; M.
\end{equation*}

Based on well-known results for Liouville-type equations (\cite{BT-2,li-cmp}), the following convergence holds
\begin{equation*}
 u_k-\bar{u}_k \ \to \sum_{j=1}^m 8\pi(1+\alpha_j)G(x,p_j) \quad {\rm in} \ \; {\rm C}_{\rm loc}^2(M\backslash \{p_1,\cdots,p_m\})
\end{equation*}
where $\bar{u}_k$ denotes the average of $u_k$ on $M$: $\overline{u}_k=\int_{M}u_k{\rm d}\mu$. For later convenience, we fix a sufficiently small $r_0>0$ and choose sets $M_j\subset M, 1\leq j\leq m$ such that
\begin{equation*}
M=\bigcup_{j=1}^m \overline{M}_j;\quad M_j\cap M_l=\varnothing,\  {\rm if}\ j\neq l;\quad B(p_j,3r_0)\subset M_j, \ j=1,\cdots,m.
\end{equation*}
Note that when $m=1$, we simply have $M_1=M$.

The concept of local maxima plays a crucial role in our analysis. For regular blow-up points (where $\alpha_j=0$), we define the maximum point $p_{k,j}$ and its corresponding value $\lambda_{k,j}$  as 
\begin{equation*}
	\lambda_{k,j}=u_k(p_{k,j}) \mathrel{\mathop:}=\max_{B(p_j,r_0)}u_k.
\end{equation*}
For singular blow-up points (where $\alpha_j\neq 0$), these quantities simplify to
\begin{equation*}
	p_{k,j}:=p_j\quad \mbox{ and }\quad \lambda_{k,j}:=u_k(p_{k,j}).
\end{equation*}
The standard bubble solution $U_{k,j}$ is defined as the radial solution to
\begin{equation}\nonumber
\Delta U_{k,j}+\rho_kh_j(p_{k,j})|x-p_{k,j}|^{2\alpha_j}e^{U_{k,j}}=0 \quad  {\rm in} \ \; \mathbb{R}^2
\end{equation}
which has the explicit form (\cite{CL1,CL2,Parjapat-Tarantello}),
\begin{equation}\nonumber
U_{k,j}(x)=\lambda_{k,j}-2\log\Big(1+\frac{\rho_k h_j(p_{k,j})}{8(1+\alpha_j)^2}e^{\lambda_{k,j}}|x-p_{k,j}|^{2(1+\alpha_j)}\Big).
\end{equation}

It is well-known (\cite{BCLT,BT-2,li-cmp}) that $u_k$ is well-approximated by the standard bubbles $U_{k,j}$ near each $p_j$: up to a uniformly bounded error term:
\begin{equation*}
  \parallel u_k(x)-U_{k,j}(x)\parallel_{L^{\infty}(B(p_j,r_0))} \leq C, 
\end{equation*}
implying, in particular, the uniform bound
\begin{equation}\nonumber
|\lambda_{k,i}-\lambda_{k,j}|\leq C, \quad 1\leq i,j \leq m.
\end{equation}

\medskip
When $\tau<m$, it has been shown in Theorem 3.3 of \cite{byz-1} that,
\begin{equation}\label{first-deriv-est}
\nabla(\log h+G_j^*)(p_{k,j})+Ce^{-\lambda_{k,j}}\phi_j^k(p_{k,j})=O(e^{-(\frac 32-\delta)\lambda_{k,j}}),\quad \tau+1\leq j\leq m,
\end{equation}
where $\phi_j^k$ is a harmonic function with bounded norms independent of $k$. Under the non-degeneracy condition $\det\big(D^2f^*(p_{\tau+1},\cdot,p_m)\big)\neq0,$ a perturbation of the Hessian matrix is still invertible. Thus we immediately know that
\begin{equation}\label{p_kj-location}
|p_{k,j}-p_j|=O(e^{-\lambda_{k,j}}),\quad \tau+1\leq j\leq m.
\end{equation} 
Later, sharper estimates were obtained in \cite{chen-lin,zhang2} for $1\leq j \leq \tau$ and
in \cite{chen-lin-sharp,gluck,zhang1} for $\tau+1\leq j\leq m$.

\medskip

\section{Local asymptotic analysis near a blow-up point}

In this section we introduce and analyze the asymptotic expansion of a sequence of blow-up solutions near both a regular point and a negative singular source. Through an appropriate coordinate transformation (with slight notational simplification), we consider a sequence of bubbling solutions $u_k$ satisfying
\begin{equation*}
\Delta u_k+|x|^{2\beta}h_k(x)e^{u_k}=0 \quad \mbox{in}\quad B_{2\tau}
\end{equation*}  
where $B_{\tau}\subset \mathbb{R}^2$ for $\tau>0$ is the ball centered at the origin with radius $\tau>0$. Here, the coefficient functions $h_k$ satisfy uniform estimates:
\begin{equation}\label{assum-H}
\frac{1}{C}\le h_k(x)\le C, \quad \|D^{m}h_k\|_{L^{\infty}(B_{2\tau})}\le C
\quad \forall x\in B_{2\tau}, \quad \forall m\le 5,
\end{equation}
and the energy have uniform bound
\begin{equation}\nonumber
\int_{B_{2\tau}}|x|^{2\beta}h_k(x)e^{u_k}<C.
\end{equation}
The sequence $\{u_k\}$ has its only blow-up point at the origin:
\begin{equation}\label{only-bu}
\mbox{there exists}\quad  z_k\to 0, \; \mbox{ such that} \quad \lim_{k\to \infty} u_k(z_k)\to \infty
\end{equation}
and for any fixed $K\Subset B_{2\tau}\setminus \{0\}$, there exists $C(K)>0$ such that
\begin{equation}\label{only-bu-2}
u_k(x)\le C(K),\quad x\in K.
\end{equation}
Additionally, it is also standard to assume that $u_k$ has bounded oscillation on $\partial B_{\tau}$:
\begin{equation}\label{u-k-boc}
|u_k(x)-u_k(y)|\le C\quad \forall x,y\in \partial B_{\tau},
\end{equation}
for some uniform  $C>0$.
To capture the boundary behavior, we introduce the harmonic function $\psi_k$ solving
\begin{equation*}
\begin{cases}
\Delta \psi_k=0  \quad \hbox{in } B_{\tau},\\
\\
\psi_k=u_k-\frac{1}{2\pi \tau}\int_{\partial B_{\tau}}u_kdS
\quad \hbox {on } \partial B_{\tau}.    
\end{cases}
\end{equation*}
The bounded oscillation condition \eqref{u-k-boc} imlpies $u_k-\psi_k$ is constant on $\partial B_{\tau}$ and $\|D^m\psi_k\|_{B_{\tau/2}}\le C(m)$ for any $m\in \mathbb N$.
The mean value property of harmonic functions also gives $\psi_k(0)=0$.

\medskip

\subsection{Asymptotic analysis around a negative source}
In this case, we use $\beta\in (-1,0)$ to denote the strength of a negative source. For  simplicity, we set
$$
 \tilde \lambda_k=u_k(0),\quad  V_k(x)=h_k(x)e^{\psi_k(x)}, \quad \tilde \varepsilon_k=e^{-\frac{u_k(0)}{2(1+\beta)}}
$$
and
\begin{equation} \label{defvk-1}
v_k(y)=u_k(\tilde \varepsilon_ky)+2(1+\beta)\log \tilde \varepsilon_k-\psi_k(\tilde \varepsilon_ky),\quad y\in \Omega_k:=B_{\tau\tilde \varepsilon_k^{-1}}(0).
\end{equation}
Obviously we have $V_k(0)=h_k(0)$ and $\Delta \log V_k(0)=\Delta \log h_k(0)$.
The standard bubble solution is given by
\begin{equation}\label{st-bub}
    U_k(y)=
  -2\log
 \left(1+\frac{h_k(0)}{8(\beta+1)^2}|y|^{2\beta+2} \right),
\end{equation}
which satisfies the equation
\begin{equation}\nonumber \Delta U_k+|y|^{2\beta}h_k(0)e^{U_k(y)}=0 \quad
 \hbox{ in } \mathbb R^2.
\end{equation}

Now we recall the following theorem:
\begin{thm}\label{integral-neg} (Theorem 3.2 of \cite{byz-1})
\begin{align*}
\int_{B_{\tau}}h_k(x)|x|^{2\beta}e^{u_k}
  =~8\pi(1+\beta)\left(1-\frac{e^{-u_k(0)}}{e^{-u_k(0)}+\frac{h_k(0)}{8(1+\beta)^2}\tau^{2\beta+2}}\right)\\
  +\tilde b_1^k(\tau,h)e^{-u_k(0)}+\tilde b_2^k(\tau,h)e^{-2u_k(0)}+O(e^{-(2+\epsilon_0)u_k(0)}).
\end{align*}
 where $\tilde b_1^k(\tau,h)$, $b_2^k(\tau,h)$ depend only on $\tau$ and derivatives of $h_k$ at $0$, and are uniformly bounded and satisfies $\lim_{\tau\to 0}\tilde b_1^k(\tau,h)=0$, $\lim_{\tau\to 0}\tilde b_2^k(\tau,h)=0$. 
\end{thm}

The main theorem in this section is on the expansion of blow-up solutions near a regular blow-up point:

\subsection{Asymptotic analysis around a regular blow-up point}
In this subsection, we analyse the solution near a regular blow-up point and consider the equation,
$$\Delta u_k+\bar h_k(x)e^{u_k}=0\quad \mbox{in}\quad B_{2\tau}. $$
We assume \[u_k(0)=\max_{x\in B_{\tau}}u_k=\bar \lambda_k,\] 
and set
\[\bar \varepsilon_k=e^{-\frac{\bar \lambda_k}{2}}.\]
We assume in particular \eqref{only-bu}, \eqref{only-bu-2},
(i.e. $0$ is the only blow-up point in $B_{2\tau}$), that the standard uniform bound
$\int_{B_{2\tau}}\bar h_ke^{u_k}\le C$ holds and finally that \eqref{u-k-boc} holds. Under these assumptions, the well-known results in \cite{li-cmp} implie that the blow-up point is simple.

Let $\phi_1^k$ be the harmonic function defined by the oscillation $u_k$ on $\partial B_{\tau}$:
\begin{equation*}
    \begin{cases}
        &\Delta \phi_1^k=0, \quad \mbox{in}\, B_{\tau},  \\
        &\phi_1^k=u_k-\frac{1}{2\pi \tau}\int_{\partial B_{\tau}}u_kdS,  \quad  \mbox{on}\, \partial B_{\tau}.
    \end{cases}
\end{equation*}
Set \[ u_1^k=u_k-\phi_1^k,\quad \mbox{in}\quad B_{\tau},\]
and 
\[h_1^k=h_ke^{\phi_1^k}\]
If we use $q_k$ to denote the maximum of $u_1^k$, it was proved in \cite{byz-1} that 
\[q_k=-\frac{2}{h_k(0)}\nabla \phi_1^k(0)\bar\varepsilon_k^2+O(\bar\varepsilon_k^3).\]
In other words,
\begin{equation}\label{qk:gluck} 
    q_k=\left(-\frac{2}{\bar h_k(0)}\partial_1\phi_1^k(0)\bar\varepsilon_k^2,-\frac{2}{\bar h_k(0)}\partial_2\phi_1^k(0)\bar\varepsilon_k^2\right)+O(\bar\varepsilon_k^3).
\end{equation}
Based on \eqref{first-deriv-est} and \eqref{qk:gluck}, we obtain the following refined estimates.

\begin{thm}\label{reg-new-e}  
\begin{equation}\label{total-q-1}
    \int_{B_{\tau}}\bar h_k(x)e^{u_k}dx=8\pi-\frac{8\pi\bar \varepsilon_k^2}{\bar\varepsilon_k^2+a_k\tau^2}-\frac{\pi\bar \varepsilon_k^2}{2}\Delta (\log \bar h_k)(0)\bar h_k(0)\bar b_{0,k}+\bar b_1^k(\tau,h)\bar \varepsilon_k^4
\end{equation}
where
\[a_k=\frac{\bar h_k(0)}{8},\]
\begin{equation*}
    \bar b_{0,k}=\int_0^{\tau\bar\varepsilon_k^{-1}}\frac{r^3(1-a_kr^2)}{(1+a_kr^2)^3}dr=\frac{1}{a_k^2}\left(1-\frac{2}{a_k\tau^2}\bar\varepsilon_k^2-\frac 12\log (a_k\tau^2)-\log \frac{1}{\bar \varepsilon_k}\right).
\end{equation*}
Here, $\bar b_1^k(\tau,h)$ depends on $\bar h_k$, $\tau$ and is uniformly bounded satisfying $\lim_{\tau\to 0+}|\bar b_1^k(\tau,h)|=0$ as $\tau\to 0$. Moreover,
for $\epsilon_0>0$ small,
\begin{equation}\label{new-rate-1}
    \nabla \log \bar h_k(q_k)+\left(1+\frac{4}{\tau\bar h_k(0)}\bar \varepsilon_k^2\right)\nabla \phi_1^k(q_k)=O(\bar\varepsilon_k^{3+\epsilon_0}),
\end{equation}
\begin{equation}\label{new-rate-2}
    \nabla \Delta \log \bar h_k(q_k)=O(\bar\varepsilon_k^{1+\epsilon_0}).
\end{equation}

\end{thm}

\begin{rem} Theorem \ref{reg-new-e} improves Theorem 3.3 of \cite{byz-1}. Even though we use the same idea of the proof,  this improvement is crucial for us to simplify the proof of the main theorem.
\end{rem}

\noindent{\bf Proof of Theorem \ref{reg-new-e}:}
We begin by employing the harmonic function $\phi_2^k$ to eliminate the oscillation of $u_1^k$ on $\partial B(q_k,\tau-|q_k|)$. According to Lemma 6.1 of \cite{byz-1}, the leading term in the expansion  of $\phi_2^k$ takes the form
$\displaystyle{-2\frac{q_k}{\tau}re^{i\theta}}$, while higher-order terms are negligible for our analysis.
 
Now we set 
\[  u_2^k=u_1^k-\phi_2^k, \quad h_2^k=h_1^ke^{\phi_2^k}.\]
We denote $\bar q_k$ to be the maximum point of $u_2^k$. Then the rough computation of $\bar q_k$ is 
\[|q_k-\bar q_k|=O(\bar\varepsilon_k^4).\]
For $u_2^k$, we mainly focus on the ball centered at $\bar q_k$ and tangent from inside to $B(q_k,\tau-|q_k|)$. Clearly the ball centered at $\bar q_k$ has a radius $\tau+O(\bar\varepsilon_k^2)$. We call this ball $B(\bar q_k,\tau_k)$. The oscillation of $u_2^k$ on $\partial B(\bar q_k,\tau_k)$ is $O(\bar\varepsilon_k^4)$.
Thus if we use $\phi_3^k$ to eliminate the oscillation of $u_2^k$ on $B(\bar q_k,\tau_k)$, the estimate of $\phi_3^k$ is
\[|\phi_3^k(\bar\varepsilon_ky)|\le C\bar\varepsilon_k^5|y|.\]

Let 
\[\bar v_k(y)=u_2^k(\bar q_k+\bar\varepsilon_ky)+2\log \bar\varepsilon_k.\]
The equation of $\bar v_k$ is
\begin{equation}\label{eq-vk-r}
    \Delta \bar v_k+ h_2^k(\bar q_k+\bar \varepsilon_ky)e^{\bar v_k}=0,\quad y\in \Omega_{2,k}:=B(0,\bar\varepsilon_k^{-1}\tau_k).
\end{equation}
By the rates obtained in \cite{byz-1}, we know
\[\nabla h_2^k(\bar q_k)=O(\bar \varepsilon_k^{3-\epsilon_0}).\]
Based on this and the approximation of $\bar v_k$ by $U_k$ we can write the equation of $\bar v_k$ as
\[\Delta \bar v_k+h_2^k(\bar q_k)e^{\bar v_k}=\left(h_2^k(\bar q_k)- h_2^k(\bar q_k+\bar\varepsilon_ky)\right)e^{\bar U_k}+O(\bar\varepsilon_k^{4+\epsilon_0})(1+|y|)^{-2-\epsilon_0},\quad \mbox{in}\quad \Omega_{2,k}.\]

We are now in a position to obtain better vanishing rates of $|\nabla h_2^k(\bar q_k)|$ and $|\nabla \Delta h_2^k(\bar q_k)|$. Let  $\Omega_{\sigma,k}=B(0,\bar\varepsilon_k^{-\sigma})$ for some $\sigma\in (0,\frac 12)$. Then we consider the following Pohozaev identity on $\Omega_{\sigma,k}$,
for any fixed $\xi\in \mathbb S^1$:
\begin{equation}\label{pi-sigma}
\begin{aligned}
    \bar\varepsilon_k\int_{\Omega_{\sigma,k}}\partial_{\xi} h_2^k(\bar q_k+\bar\varepsilon_ky)e^{\bar v_k}=&\int_{\partial \Omega_{\sigma,k}}\left(h_2^k(\bar q_k+\bar\varepsilon_ky)e^{\bar v_k}(\xi,v)+\partial_{v}\bar v_k\partial_{\xi}\bar v_k\right)\\
    &-\frac 12\int_{\partial \Omega_{\sigma,k}} |\nabla \bar v_k|^2(\xi\cdot v). 
\end{aligned}
\end{equation}
To evaluate the left hand side of \eqref{pi-sigma}, we use the expansion of $\bar v_k$ and $\partial_{\xi} h_2^k(x)$.  Note that if we use $H_{\xi}^k$ to denote all the high order terms in the expansion of $\partial_{\xi} h_2^k$:
\begin{equation*}
    \partial_{\xi} h_2^k(x)=\partial_{\xi}h_2^k(\bar q_k)+\sum_{j=1}^3\sum_{|\alpha |=j}\frac{1}{\alpha !}\partial_{\xi}\partial^{\alpha} h_2^k(\bar q_k)(x-\bar q_k)^{\alpha}+H_{\xi}^k(x),
\end{equation*}
we can choose $\sigma$ close to $\frac 12$ so the integration with $H_{\xi}^k$ is minor compared to the error $O(\bar\varepsilon_k^{4+\epsilon_0})$.  Then after using the current vanishing rate of $|\nabla  h_2^k(0)|=O(\bar\varepsilon_k^{3-\epsilon_0})$, and taking $\sigma$ close to $\frac 12$ but less than $\frac 12$, we have 
\begin{equation*}
    LHS=8\pi\bar\varepsilon_k\partial_{\xi}\log h_2^k(\bar q_k)+\frac{16\pi}{ h_2^k(\bar q_k)^2}\partial_{\xi}\Delta h_2^k(\bar q_k)\bar\varepsilon_k^3(\log \bar\varepsilon_k^{-2\sigma}+C_k)+O(\bar\varepsilon_k^{4+\epsilon_0}),
\end{equation*}
with some bounded sequence \(C_k\).

To estimate the right hand side of \eqref{pi-sigma}, by using the expansion of $\bar v_k$ and $h_2^k$ with $\sigma$ close to $\frac 12$, the first term $\int_{\partial \Omega_{\sigma,k}} h_2^ke^{v_k}(\xi\cdot v)$ gives
\[\int_{\partial \Omega_{\sigma,k}} h_2^ke^{\bar v_k}(\xi\cdot v)=C_{1,k}\bar\varepsilon_k^3\partial_{\xi}\Delta h_2^k(\bar q_k)+O(\bar\varepsilon_k^{4+\epsilon_0}),\]
with some bounded sequence $C_{1,k}$, while in the last two boundary terms, because of symmetry and cancellations, the only term that needs to be evaluated is 
 \begin{equation}\label{est-gra}
     \int_{\partial \Omega_{\sigma,k}}\bar U_k'(r)\partial_{\xi}c_{1,k}.
 \end{equation}
From Theorem 3.3 and its proof, we know that
\[
  |c_{1,k}(y)|\leq C \varepsilon_k^{4-\epsilon_0}\frac{r^2}{1+r}
\]
We readily have the following fact
\[\bar U_k'(r)=-\frac{4 h_2^k(0)r}{1+ h_2^k(0)r^2},\]
\[\partial_{\xi}c_{1,k}=O(\bar\varepsilon_k^{4-\epsilon_0})\quad \text{on} \; \partial \Omega_{\sigma,k},\]
This means the leading term of the integral \eqref{est-gra} is $E_k\bar\varepsilon_k^{4-\epsilon_0}$ with some bounded sequence $E_k$ independent of $\sigma$. The error term is $O(\bar\varepsilon_k^{4+2\sigma-\epsilon_0})$. Thus, the choice of $\sigma$ insures the error term of the integral is $O(\bar\varepsilon_k^{4+\epsilon_0})$. At this point, we have
\[RHS=C_1\bar\varepsilon_k^3\partial_{\xi}\Delta  h_2^k(\bar q_k)+C_2\bar\varepsilon_k\partial_{\xi}h_2^k(\bar q_k)+E_k\bar\varepsilon_k^{4-\epsilon_0}+O(\bar\varepsilon_k^{4+\epsilon_0}).\] 
Finally, the combination of both sides gives
\begin{equation}\label{rate-1}
        8\pi\bar\varepsilon_k\partial_{\xi}\log h_2^k(\bar q_k)+\frac{16\pi}{ h_2^k(\bar q_k)^2}\bar\varepsilon_k^3\partial_{\xi}\Delta h_2^k(\bar q_k)(\log \bar\varepsilon_k^{-2\sigma}+C_k)+E_k\bar\varepsilon_k^{4-\epsilon_0}=O(\bar\varepsilon_k^{4+\epsilon_0}).
\end{equation} 
Therefore by choosing $\sigma=\sigma_1$ and $\sigma_2$ with both of them close to and greater than $\frac 12$ (without affecting the error term $O(\bar\varepsilon_k^{4+\epsilon_0})$, we have
\[|(\sigma_1-\sigma_2)\partial_{\xi}\Delta h_2^k(\bar q_k)|=O(\bar\varepsilon_k^{1+\epsilon_0})\]
for $\epsilon_0>0$ arbitrarily small.
Since $\xi$ is any vector on $\mathbb S^1$, we have
\begin{equation}\label{rate-2}
\nabla\Delta  h_2^k(\bar q_k)=O(\bar\varepsilon_k^{1+\epsilon_0}),
\end{equation}
which immediately implies that 
\begin{equation}\label{rate-3}\nabla  h_2^k(\bar q_k)=O(\bar\varepsilon_k^{3+\epsilon_0}),
\end{equation}
by the Pohozaev identity \eqref{rate-1}. 
Tracking the definition of $h_2^k$, we see that (\ref{rate-3}) is equivalent to 
\[\nabla \log \bar h_k(q_k)+\nabla \phi_1^k(q_k)(1+\frac{4}{\tau \bar h_k(q_k)}\bar\varepsilon_k^2)=O(\bar\varepsilon_k^{3+\epsilon_0}).\]
Now the expansion of $v_k$ can be written as
\begin{equation}\label{vk-expan}
\bar v_k(y)=U_k(y)+\sum_{l=0}^4c_{l,k}+O(\bar\varepsilon_k^{4+\epsilon_0})(1+|y|)^{2+\epsilon_0},\quad |y|\le \tau \bar\varepsilon_k^{-1}.
\end{equation}
By incorporating the improved decay rates of both first and third derivatives in our integral estimates, we establish precise bounds for $b_0^k$ and $b_1^k$. Finally, we remark that even though we did not evaluate the integral on a ball centered at the origin, the behavior of $v_k$ outside the symmetry area is of the magnitude $O(\bar\varepsilon_k^4)$ and the difference between the two balls is $O(\bar\varepsilon_k)$. These boundary effects are negligible (higher-order) in the final estimates
\qed


\section{Proof of the Main Theorem}

In this section, we prove Theorem \ref{mainly-case-2}.
We establish the result by contradiction, considering two distinct blow-up sequences $v_1^k$ and $v_2^k$ that solve equation (\ref{m-equ}) corresponding to the same $\rho_k$ and the same set of blow-up points. The case we consider in this context is $\alpha_M\le 0$.
For clarity, we initially examine the case featuring exactly two blow-up points: a regular point $q$ and a negative singular source $p$ with the strength $\beta\in (-1,0)$. This representative case captures all essential analytical challenges while allowing straightforward generalization to configurations involving multiple regular points and additional negative singular sources through analogous arguments. The proof strategy relies on careful comparison of the asymptotic behaviors of these competing blow-up sequences near the distinguished points $p$ and $q$.

Let
$$w_i^k=v_i^k+4\pi \beta G(x,p),\quad i=1,2.$$
Then we have
\begin{equation}\label{uik-tem}
\Delta_g w_i^k+\rho_k\left(\frac{He^{w_i^k}}{\int_MHe^{w_i^k}}-1\right)=0,
\end{equation}
where
$$H(x)=h(x) e^{-4\pi\beta G(x,p)}.$$
Since by adding any constant to $w_i^k$ in (\ref{uik-tem}) we still come up with a solution, we fix this free constant to make
$w_i^k$ satisfy
\begin{equation*}
\int_M He^{w_i^k}dV_g=1,\quad i=1,2.
\end{equation*}
The basic assumption is that
$$\sigma_k:=\|w_1^k-w_2^k\|_{L^{\infty}(M)}>0.$$
From previous works such as \cite{bart-4}, \cite{wu-zhang-ccm} we have that $\rho_1^k=\rho_2^k=\rho_k$
implies that $\sigma_k=o(1)$ where we set $\hat \lambda_i^k=w_i^k(p)$ and $\bar \lambda^k_i=\max\limits_{B(q,\tau)}w_i^k$, it is well known that
\begin{equation}\label{same-ord}
    \hat \lambda^k_i-\bar \lambda^k_i=O(1),\quad i=1,2.\quad |\tilde \lambda_1^k-\tilde \lambda_2^k|\le C(\tilde \lambda_1^k)^{-1},\quad |\bar \lambda_1^k-\bar \lambda_2^k|\le C(\bar \lambda_1^k)^{-1}.
\end{equation}

We work in local isothermal coordinates around $p,q$, the metric locally near $p$ taking the form,
$$ds^2=e^{\phi_p}((dx_1)^2+(dx_2)^2),$$ where $\phi_p$ satisfies,
\begin{equation}\label{def-phi}
\Delta \phi_p+2Ke^{\phi_p}=0 \quad{\rm in }\quad B_{\tau}, \quad \phi_p(0)=|\nabla \phi_p(0)|=0.
\end{equation}
In the same fashion, we can define $\phi_q$.
Remark that in our setting, we have,
\begin{equation}\label{rho-3}
\begin{cases}
\rho^*=16\pi+8\pi\beta,\quad N^*=4\pi\beta,\\
\\
G^*_p(x)=8\pi(1+\beta)R(x,p)+8\pi G(x,q). 
\end{cases}
\end{equation}
Similarly, we set $G_q^*$ as
\[G_q^*(x)=8\pi(1+\beta)G(x,p)+8\pi R(x,q).\]
\begin{rem}
We can always choose $\tau$ small enough to guarantee that
$\Omega(p,\tau)$ and $\Omega(q,\tau)$ are simply connected and at
positive distance one from each other. After scaling,
$x=\varepsilon y$ for some $\varepsilon>0$, we will denote
$B_{\tau/\varepsilon}=B_{\tau/\varepsilon}(0)$. Moreover, for notational simplicity, we intentionally reuse function symbols $w_k$, $\xi_k$, $\tilde h_k$ across different local coordinate systems. Of course, it will be clear, time to time, by the context what the given symbol means.
\end{rem}
\noindent
Thus, working in these local coordinates centered at $p$, $q$ respectively, we have the local variables $x\in B_\tau$ around $p$ and we define $f_{p,k}$ to be the solution of,
\begin{equation}\label{def-fk}
\Delta f_{p,k}=e^{\phi_p}\rho_k\quad {{\rm in}}\quad B_{\tau}, \quad f_{p,k}(0)=0.
\end{equation}
$f_{q,k}$ is defined similarly.
To write the local equation around 
 $p$, we set
\begin{equation}\label{uk-def1}
\tilde u_i^k= w_i^k-f_{p,k} \quad { {\rm in} }\quad B_\tau.
\end{equation}
Then we have 

\begin{equation}
\label{around-p2-h}
\begin{cases}
\Delta \tilde u_i^k+|x|^{2\beta}h_p^k e^{\tilde u_i^k}=0\quad {{\rm in}}\quad B_{\tau},\quad i=1,2,\qquad 0=x(p)\\
\\
h_p^k=\rho_k h(x)e^{-4\pi \beta R(x,p)+\phi_p+f_{p,k}}.
\end{cases}
\end{equation}
Then around $q$ we set
\[\bar u_i^k=w_i^k-f_{q,k}.\]
Then the equation for $\bar u_i^k$ is
\begin{equation}
\label{around-q}
\begin{cases}
\Delta \bar u_i^k+h_q^k e^{\bar u_i^k}=0\quad {{\rm in}}\quad B_{\tau},\quad i=1,2,\qquad 0=x(q),\\
\\
h_q^k=\rho_k h(x)e^{-4\pi \beta G(x,p)+\phi_q+f_{q,k}}.
\end{cases}
\end{equation}
Obviously $\hat h_k, \bar h_k$ satisfy \eqref{assum-H}.

\begin{rem}\label{center-q}
For technical convenience in handling \eqref{around-q}, it will be useful to define $\bar q_i^k\to 0$ to be the maximum points of $u_i^k$ in $B_\tau$ and work 
in the local coordinate system centered at $\bar q_i^k$. In particular, $q_i^k\in M$ will denote the pre-images of these points via the local isothermal map, and we will sometimes work, possibly taking a smaller $\tau$, with \eqref{around-q} where $0=x(q^k_i)$. 
\end{rem}
By using also \eqref{def-phi}, \eqref{def-fk}, we see that the non-degeneracy assumption $L(\mathbf{p})\neq 0$ takes the form
$$\Delta_g \log  h(q)+\rho_*-N^*-2K(q)\neq 0, $$
which, for any $k$ large enough, is equivalent to (see (\ref{rho-3}))
$$\Delta \log h_q^k(0)\neq 0.$$ 
Also, because of the new vanishing rate of the first derivatives (\ref{new-rate-1}), \eqref{same-ord} and the non-degeneracy assumption $det(D^2f^*(p_{\tau+1},...,p_m))\neq 0$, we see that we have a good estimate on $|q_1^k-q_2^k|$:
\begin{equation}\label{close-q}
|q_1^k-q_2^k|= o(e^{-\bar\lambda_i^k} )
\end{equation}
because the Hessian matrix is invertible, and it is still invertible under a perturbation, which comes from the extra term in (\ref{new-rate-1}). 

At this point, let us define,
\begin{equation*}
\xi_k:= \frac{w_1^k-w_2^k}{\sigma_k}.
\end{equation*}
Remark that, because of \eqref{uk-def1}, locally we have
\begin{equation}\label{ueqw}
\tilde u_1^k-\tilde u_2^k= w_1^k-w_2^k,\quad \bar u_1^k-\bar u_2^k= w_1^k-w_2^k
\end{equation}
and in particular in local coordinates $\xi_k=(u_1^k-u_2^k)/\sigma_k$. Then it is well known (see \cite{bart-4,bart-4-2,wu-zhang-ccm}) that, after a suitable scaling in local coordinates, the limit of $\xi_k$ takes the form,
\begin{align}
    &b_0\frac{1-B|y|^{2\beta+2}}{1+B|y|^{2\beta+2}}  &&  \mbox{\rm near } p, \quad B=\lim_{k\to \infty}\frac{h_p^k(0)}{8(1+\beta)^2};  \label{linim-2}  \\
    &b_0\frac{1-C|y|^2}{1+C|y|^2}+b_1\frac{y_1}{1+C|y|^2}+b_2\frac{y_2}{1+C|y|^2}  &&  \mbox{\rm near } q,   \quad C=\lim_{k\to \infty}\frac{ h_q^k(0)}8.   \label{linim-3}
\end{align}

\begin{rem}\label{xi:b_0}
The reason why we find exactly the same $b_0$ in \eqref{linim-2} and \eqref{linim-3} is that $\xi_k$ converges to a constant far away from blow-up points,
and it is by now well known (see \cite{bart-4,bart-4-2,wu-zhang-ccm})
that this constant is exactly $-b_0$.
\end{rem}

\subsection{ First estimates about $\sigma_k$. }

As mentioned above, the maximum of $\tilde u_i^k$ around $p$ is $\tilde \lambda_i^k$ ($\tilde \lambda_{i,k}=\tilde u_i^k(p)$ ). The maximum of $\bar u_i^k$ around $q$ is $\bar \lambda_i^k$ ($\bar \lambda_i^k=\bar u_i(q_i^k)$).
We derive first a rough estimate for $w_1^k$ and $w_2^k$ far away from blow-up points:

\begin{lem}\label{w12-away} 
For $x\in M\setminus (\Omega(p,\tau)\cup \Omega(q,\tau))$,
\begin{equation}\label{uk1}
\begin{aligned}
   w_i^k(x)=& \bar{w}_i^k+8\pi(1+\beta)(1-\frac{e^{-\tilde \lambda_i^k}}{e^{-\tilde \lambda_i^k}+\frac{h_p^k(0)}{8(1+\beta)^2}\tau^{2\beta+2}})G(x,p)\\
   & +\bigg(8\pi(1-\frac{e^{-\bar \lambda_i^k}}{e^{-\bar \lambda_i^k}+a_k\tau^{-2}})-\frac{\pi}2\Delta (\log h_q^k(0))h_q^k(0)\bar b_0^k e^{-\bar \lambda_i^k}\bigg )G(x,q_i^k)\\
   &+c_k(\tau,h)e^{-\bar \lambda_i^k}+O(e^{-(1+\epsilon_0)\bar\lambda_i^k}),\quad i=1,2,
\end{aligned}
\end{equation}
where $q_i^k$ is defined in Remark \ref{center-q},$\epsilon_0>0$ is a small constant, $|c_k(\tau,h)|\to 0$ as $\tau\to 0$.
\end{lem}

\noindent{\bf Proof of Lemma \ref{w12-away}:}
From the Green representation formula (recall that $\int_M G(x,\eta)=0$) we have,
\begin{align*}
w_i^k(x)-\bar{w}_i^k&=\int_M G(x,\eta)\rho_kH e^{w_i^k}d\mu\\
&=G(x,p)\int_{\Omega({p},\tau)}\rho_kH e^{w_i^k}d\mu+
\int_{\Omega({p},\tau)}(G(x,\eta)-G(x,p))\rho_kH e^{w_i^k}  \\
&\quad+G(x,q_i^k)\int_{\Omega(q_1^k,\tau)}\rho_kH e^{w_i^k}
+\int_{\Omega(q_i^k,\tau)}(G(x,\eta)-G(x,q_i^k))\rho_kH e^{w_i^k}\\
&\quad+\int_{M\setminus \{ \Omega({p},\tau)\cup
\Omega(q_i^k,\tau)\}}G(x,\eta)\rho_kH e^{w_i^k}.
\end{align*}

By using the expansion theorems (Theorem \ref{integral-neg} and Theorem \ref{reg-new-e}), we have 
\begin{align*}
&\int_{\Omega({p},\tau)}\rho_kH e^{w_i^k}d\mu=\int_{B_{\tau}} |x|^{2\beta}h_p^k e^{\tilde u_i^k}dx\\
=&8\pi(1+\beta)(1-\frac{e^{-\tilde \lambda_{i,k}}}{e^{-\tilde \lambda_{i,k}}+\frac{h_p^k(0)}{8(1+\beta )^2}\tau^{2\beta+2}})+\tilde b_1^k(\tau,h)e^{-\tilde \lambda_i^k}+\tilde b_2^k(\tau,h)e^{-2\tilde \lambda_i^k}+O(\bar\varepsilon_k^{4+\epsilon_0}),
\end{align*}
and
\begin{align*}
&\int_{\Omega(q_i^k,\tau)}\rho_kH e^{w_i^k}
=\int_{B_{\tau}} h_q^k e^{\bar u_i^k}dx\\
=&8\pi-\frac{8\pi e^{-\bar \lambda_i^k}}{e^{-\bar \lambda_i^k}+\frac{\bar h_q^k(0)}8\tau^2}-\frac{\pi e^{-\bar \lambda_i^k}}{2}\Delta (\log h_q^k(0)h_q^k(0)\bar b_0^k+\bar b_1^k(\tau,h) e^{-\bar \lambda_i^k}
\end{align*}
where $\bar b_0^k$ is defined in Theorem \ref{reg-new-e}, $\lim_{\tau\to 0}\bar b_1^k(\tau,h)=0$. 

We sketch the argument behind the
estimates of the second and fourth terms.
Using the expansion of blow-up solutions, we obtain
\[\int_{\Omega({p},\tau/2)} (G(x,\eta)-G(x,p))\rho_kH e^{w_i^k}=\tilde b_3^k(h,\tau)e^{-\tilde \lambda_1^k}+O(e^{-(1+\epsilon_0)\tilde \lambda_i^k}),\]
where $\tilde b_3^k(\tau,h)\to 0$ as $\tau \to 0$. This computation is carried out in local coordinates $y(\eta)$;
we use the Taylor expansion of $G(x,\eta_y)$ to $y$ around $0=y(p_i)$ and the
expansion of $u^k_1$. Then it is easy to see that the
$O(\varepsilon_k)$ term cancels out due to the integration on the symmetric set $B_{\tau}$
and the leading term is of the order 
$O(\varepsilon_k^2)$. In the integration of $\int_{\Omega(q_1^k,\tau/2)}(G(x,\eta)-G(x,q_i^k))\rho_kHe^{w_i^k}$,
we adopt the same argument with minor changes, this time based on Theorem 1.2 in \cite{zhang1},
which yields as well to an error of order $O(\varepsilon_k^{2+\epsilon_0})$ for some $\epsilon_0>0$ easily.
This fact concludes the $w_i$ part of the proof of \eqref{uk1}. 

In the last integral, using the fact $\bar w_i^k=-\tilde\lambda_i^k+O(1)$, we obtain the rough estimate 
\[\int_{M\setminus \{ \Omega({p},\tau)\cup
\Omega(q_i^k,\tau)\}}G(x,\eta)\rho_kH e^{w_i^k}=O(e^{-\bar\lambda_i^k}).\]
Therefore, the above analysis yields
\[w_i^k(x)=\bar{\rm{w}}_i^k-4(1+\beta)\log |x-p|+G_p^*(x)+O(e^{-\tilde \lambda_1^k}).\]
On the other hand, using the local coordinates around $p$,  we have
\begin{equation*}
    \begin{aligned}
\tilde u_i^k(x)=-\tilde \lambda_i^k-4(1+\beta)\log |x|+O(1)+O(e^{-\tilde \lambda_1^k}).  
\end{aligned}
\end{equation*}
Comparing the two estimates above, we have that 
\begin{equation*}
\bar w_i^k=-\tilde \lambda_i^k-2\log \frac{ h_p^k(0)}{8(1+\beta)^2}-G_p^*(p)+O(e^{-\lambda_1^k}).
\end{equation*}
Hence, we further derive\[\int_{M\setminus \{ \Omega({p},\tau)\cup
\Omega(q_i^k,\tau)\}}G(x,\eta)\rho_kH e^{w_i^k}=O(e^{-(1+\epsilon_0)\bar\lambda_i^k}).\]
By summing up the estimates of the above five terms, we complete the proof of Lemma \ref{w12-away}. 

\qed

\medskip

Here and in the rest of the proof, we denote by,
$$c_k(x)=\begin{cases}
\dfrac{e^{w_1^k(x)}-e^{w_2^k(x)}}{w_1^k(x)-w_2^k(x)},\quad &\mbox{if} \quad w_1^k(x)\neq w_2^k(x),\\
\\
e^{w_1^k(x)}, \quad &\mbox{if }\quad w_1^k(x)=w_2^k(x),
\end{cases}
$$
and in local coordinates near $p$ we define $\tilde c_k$ as 
$$\tilde c_k(x)=\begin{cases}
\dfrac{e^{u_1^k(x)}-e^{u_2^k(x)}}{u_1^k(x)-u_2^k(x)},\quad &\mbox{if} \quad u_1^k(x)\neq u_2^k(x),\\
\\
e^{u_1^k(x)}, \quad &\mbox{if }\quad u_1^k(x)=u_2^k(x),
\end{cases}
$$
where (see \eqref{def-fk}),
\begin{equation}\label{ceqtc}
c_k(x)=e^{f_p^k(x)}\tilde c_k(x).
\end{equation}
We can define $\bar c_k$ near $q$ similarly. Near $q$ we have 
$c_k(x)=e^{f_q^k}\bar c_k$. 
Now we prove a point-wise estimate up to the order $o(e^{-\bar \lambda_i^k})$.
\begin{lem}\label{point-compare}
\begin{align}
    -\bar \lambda_i^k&=-\tilde \lambda_i^k+2\log\frac{h_q^k(0)}{8}-2\log \frac{h_p^k(0)}{8(1+\beta)^2}+G_q^*(q)-G_p^*(p)+ b_w^k(\tau,h)e^{-\tilde \lambda_1^k}+O(e^{-(1+\epsilon_0)\lambda_1^k}).  \label{difference-lambda}  \\
    e^{-\bar \lambda_i^k}&=e^{-\tilde\lambda_i^k}\frac{(1+\beta)^4(h_q^k)^2(0)}{(h_p^k)^2(0)}e^{G_q^*(q)-G_p^*(p)}+ b_w^k(\tau,h)e^{-2\tilde \lambda_1^k}+O(e^{-(2+\epsilon_0)\lambda_1^k}).  \label{height-d}
\end{align}
Here, $b_w^k(\tau,h)$ may denote the different bounded functions related to $\tau$ and $h$.

\end{lem}

\noindent{\bf Proof of Lemma \ref{point-compare}:}Here we don't specify $q_1^k$ or $q_2^k$ because $|q_1^k-q_2^k|=o(\bar\varepsilon_k^2)$.  We carry out the analysis around $p$. 

Let $\phi_{i,p,k}$ be the harmonic function that eliminates the oscillation of $\tilde u_i^k$ on $\partial B(p,\tau)$ and take $p$ as the origin, then $\tilde u_i^k$ in $B(p,\tau)$ is
\begin{align*}
&(\tilde u_i^k)(p+\tilde \varepsilon_i^ky)+2(1+\beta)\log \tilde \varepsilon_i^k\\
=&U_i^k(\tilde \varepsilon_i^ky)+\phi_{i,p,k}(\tilde \varepsilon_i^ky)+\sum_{s=0}^2c_{s,k}+O(\tilde \varepsilon_k^{3+3\beta+\epsilon_0})(1+|y|)^{3+3\beta+\epsilon_0},\quad x\in B_{\tilde \varepsilon_k^{-1}\tau}
\end{align*}
where $c_{s,k}$ is the projection of the function on the left to $e^{is\theta}$,
\[U_i^k(x)=\tilde \lambda_i^k-2\log\left(1+\frac{h_p^k(0)}{8(1+\beta)^2}e^{\tilde \lambda_i^k}|x|^{2+2\beta}\right),\]
$\tilde \varepsilon_i^k=e^{-\frac{\tilde \lambda_i^k}{2(1+\beta)}}$. In particular for $x\in B_{\tau}\setminus B_{\tau/2}$, we have
\begin{equation}\label{uk-neg-d}
\begin{aligned}
\tilde u_i^k(x)=&-\tilde \lambda_i^k-2\log \frac{h_p^k(0)}{8(1+\beta)^2}-4(1+\beta)\log |x|\\
&+\phi_{i,p,k}(x)+\sum_{s=0}^2c_{s,k}+O(e^{-(1+\epsilon_0)\tilde \lambda_i^k})  
\end{aligned}
\end{equation}
for some $\epsilon_0>0$. Here for convenience we define 
\begin{equation}\label{big-phi}
\Phi(x)= 8\pi(1+\beta)G(x,p)+8\pi G(x,q).
\end{equation}
On the other hand, the Green's representation of $w_i^k$ gives
\begin{equation}\label{green-u-away}
w_i^k(x)=\bar w_i^k+\Phi(x)+b_w^k(h,\tau)e^{-\bar \lambda_i^k}+O(e^{-(1+\epsilon_0)\tilde \lambda_i^k})
\end{equation}
  for $x$ away from singular sources, where $\Phi$ is defined in (\ref{big-phi}) and $x\in B(p,2\tau)$. By the definition of $G_p^*$, this expression can be written as 
\[w_i^k(x)=\bar{\rm{w}}_i^k-4(1+\beta)\log |x-p|+G_p^*(x)+b_w^k(\tau,h)e^{-\tilde \lambda_1^k}+O(e^{-(1+\epsilon_0)\tilde \lambda_1^k}).\]
Inserting $f_{p,k}$ and using the local coordinates around $p$,  we have
\begin{equation}\label{uk-ave}
    \begin{aligned}
\tilde u_i^k(x)=&-\tilde \lambda_i^k-4(1+\beta)\log |x|+G_p^*(p)\\
&+G_p^*(x)-G_p^*(p)-f_{p,k}+ b_w^k(\tau,h)e^{-\tilde \lambda_1^k}+O(e^{-(1+\epsilon_0)\tilde \lambda_i^k}).  
\end{aligned}
\end{equation}

Comparing (\ref{uk-neg-d}) and (\ref{uk-ave}) we have that 
\begin{equation}\label{uk-bar-l}
\bar w_i^k=-\tilde \lambda_i^k-2\log \frac{ h_p^k(0)}{8(1+\beta)^2}-G_p^*(p)+ b_w^k(\tau,h)e^{-\tilde \lambda_1^k}+O(e^{-(1+\epsilon_0)\lambda_1^k}). 
\end{equation}
Around $q$ we have
\begin{equation}\label{uk-bar-1}
\bar w_i^k=-\bar \lambda_i^k-2\log \frac{h_q^k(0)}{8}-G_q^*(q)+ b_w^k(\tau,h)e^{-\tilde \lambda_1^k}+O(e^{-(1+\epsilon_0)\tilde \lambda_i^k}).
\end{equation}

Equating (\ref{uk-bar-1}) and (\ref{uk-bar-l}), we have
\[
-\bar \lambda_i^k=-\tilde \lambda_i^k+2\log\frac{h_q^k(0)}{8}-2\log \frac{h_p^k(0)}{8(1+\beta)^2}+G_q^*(q)-G_p^*(p)+ b_w^k(\tau,h)e^{-\tilde \lambda_1^k}+O(e^{-(1+\epsilon_0)\tilde \lambda_i^k}).
\]
Consequently (\ref{height-d}) holds. 
Lemma \ref{point-compare} is established. 
\qed

\medskip

The following theorem is inspired by Theorem 2.2 in \cite{byz-2}. 

\begin{thm}\label{combine-rho} For $i=1,2$,
\begin{align*} 
\rho_k-\rho_*=&\frac{16\pi e^{-\bar \lambda_i^k}}{h_p^k(0)^2e^{G_q^*(q)}}L(\mathbf{p})\left(\bar\lambda_i^k+\log (h_q^k(0)^2e^{G_q^*(q)}\tau^2)-2\right)\\
&+\frac{64 e^{-\bar \lambda_i^k}}{h_q^k(0)^2}(D(\mathbf{p})+o(1))+B_k(\tau,h)e^{-2\bar  \lambda_i^k}
+O(e^{-(2+\epsilon_0)\bar\lambda_i^k}).
\end{align*}
where 
$o(1)\to 0$ as $\tau\to 0$.
If both regular and singular blow-up points exist, $\bar \lambda_i^k$ and $\tilde \lambda_i^k$ are related by \eqref{difference-lambda} ($i=1,2$).
 \end{thm}

\noindent{\bf Proof of Theorem \ref{combine-rho}:} We evaluate 
$\rho_k=\int_M\rho_kHe^{\rm w_i^k}d\mu$.
\begin{align*}
\rho_k=\int_{B(p,\tau)}\rho_kHe^{\rm w_i^k}d\mu+\int_{B(q,\tau)}\rho_kHe^{\rm w_i^k}d\mu+\int_E\rho_kHe^{\rm w_i^k}d\mu
\end{align*}
where $E=M\setminus (B(p,\tau)\cup B(q,\tau))$. Let $\Omega_p$ be a neighborhood of $p$, $\Omega_q$ a neighborhood of $q$ such that $\Omega_p\cup \Omega_q=M$, $\Omega_p\cap \Omega_q=\emptyset$. By means of Theorem \ref{integral-neg}, we have
\begin{equation}
\label{inter-ne-1}
\begin{aligned}
&\int_{B(p,\tau)}\rho_kHe^{\rm w_i^k}d\mu
=\int_{B_{\tau}}\rho_k|x|^{2\beta}h_p^ke^{f_{p,k}+\phi_p}e^{\tilde u_i^k}dx \\
=& 8\pi(1+\beta)\left(1-\frac{e^{-\tilde \lambda_i^k}}{e^{-\tilde\lambda_i^k}+\frac{h_p^k(0)}{8(1+\beta)^2}\tau^{2+2\beta}}\right)+\tilde b_1^k(\tau,h)e^{-\tilde \lambda_i^k}+\tilde b_2^k(\tau,h)e^{-2\tilde \lambda_i^k}+O(e^{-(2+\epsilon_0)\tilde\lambda_i^k}),
\end{aligned}
\end{equation}
where $\tilde b_1^k(\tau,h)\to 0$ and $\tilde b_2^k(\tau,h)$ as $\tau\to 0$.
For the integral outside the $B(p,\tau)$, by (\ref{green-u-away}) and (\ref{uk-bar-l}),
\begin{align*}
&\int_{\Omega_p\setminus B(p,\tau)}\rho_kHe^{\rm w_i^k}d\mu\\
=&\int_{\Omega_p\setminus B(p,\tau)}\rho_kHe^{\bar w_i^k}e^{\Phi}d\mu+O(e^{-(2+\epsilon_0)\tilde \lambda_i^k})\\
=& e^{-\tilde \lambda_i^k}\left(\frac{8(1+\beta)^2}{h_p^k(0)}\right)^2e^{-G_p^*(p)}\int_{\Omega_p\setminus B(p,\tau)}\rho_kHe^{\Phi}d\mu+B_k(\tau,h)e^{-2\tilde \lambda_i^k}+O(e^{-(2+\epsilon_0)\tilde \lambda_i^k}).
\end{align*}
Here, $B_k(\tau,h)$ can denote the different bounded functions related to $\tau$ and $h$. When we use local coordinates around $p$ to evaluate the last integral, we see that 
\begin{align*}\int_{\Omega_p\setminus B(p,\tau)}\rho_kHe^{\Phi}d\mu 
=&\int_{\Omega_p\setminus B(p,\delta)}\rho_kH e^{\Phi}d\mu+\int_{B(p,\delta)\setminus B(p,\tau)}h_p^k|x|^{-4-2\beta}e^{G_p^*(x)}dx\\
=&\int_{\Omega_p\setminus B(p,\delta)}\rho_kH e^{\Phi}d\mu+
h_p^k(0)e^{G_p^*(p)}\left(\frac{\pi\tau^{-2-2\beta}}{1+\beta}+O(\delta^{-2-2\beta})\right).
\end{align*}
Therefore 
\begin{align*}
    \int_{\Omega_p}\rho_kHe^{w_i^k}d\mu
=&8\pi(1+\beta)\left(1-\frac{e^{-\tilde \lambda_i^k}}{e^{-\tilde\lambda_i^k}+\frac{h_p^k(0)}{8(1+\beta)^2}\tau^{2+2\beta}}\right) \\
&+e^{-\tilde\lambda_i^k}e^{-G_p^*(p)}\left(\frac{8(1+\beta)^2}{h_p^k(0)}\right)^2\left(D_s+\tilde b_3^k(\delta,h)\right)\\
&+\tilde b_4^k(\delta,h)e^{-2\tilde \lambda_i^k}+O(e^{-(2+\epsilon_0)\tilde\lambda_i^k}),
\end{align*}
where $|\tilde b_3^k(\delta,h)|\le C\delta^{-2-2\beta}$, 
$$D_s=\lim_{\tau\to 0+}\int_{\Omega_p \setminus B(p,\tau)} \rho_kHe^{\Phi}d\mu-\frac{h_p^k(0)e^{G_p^*(p)}\pi}{1+\beta}\tau^{-2-2\beta}.$$

Now, for integration around $q$, we need to deal with a logarithmic term in integration. For this, we set 
\[\tau_q=\sqrt{8 e^{G_q^*(q)}h_q^k(0)}\tau, \]
and invoke \eqref{total-q-1} of Theorem \ref{reg-new-e}: 
\begin{align*}
&\int_{B(q,\tau_q)}\rho_kHe^{\rm w_i^k}d\mu=\int_{B_{\tau_q}}h_q^ke^{f_{q,k}+\phi_q}e^{\bar u_i^k}dx\\
&=8\pi-8\pi \frac{e^{-\bar\lambda_i^k}}{e^{-\bar\lambda_i^k}+a_k\tau_q^2}-\frac{\pi e^{-\bar\lambda_i^k}}{2}\Delta (\log h_q)(0)h_q^k(0)\bar b_{0,k}+\bar b_1^k(\tau,h)e^{-\bar\lambda_i^k}
.
\end{align*}
where $\epsilon_0>0$ is a small constant, $\bar b_1^k(\tau,h)\to 0$ when $\tau\to 0$, and 
\[\bar b_{0,k}=\int_0^{\tau_q \bar \varepsilon_k^{-1}}\frac{r^3(1-a_kr^2)}{(1+a_kr^2)^3}dr,\quad a_k=\frac{ h_q^k(0)}8.\]
Here we use elementary method to estimate $b_{0,k}$ as
\[\bar b_{0,k}=\frac{1}{2a_k^2}(-\bar\lambda_i^k-\log (a_k\tau_q^{2})+2)+d_k(\tau,h)e^{-\bar\lambda_i^k},\]
where $d_k(\tau,h)$ depend only on $\tau$ and $h$ and satisfies $\lim_{\tau\to0}d_k(\tau,h)=0$. Using this expression of $b_{0,k}$ we have, 
\begin{equation}
\label{ener-reg-2}
\begin{aligned}
\int_{B(q,\tau_q)}\rho_kHe^{\rm w_i^k}d\mu
=&8\pi -\frac{16\pi}{h_q^k(0)}\Delta \log h_q^k(0)e^{-\bar\lambda_i^k}\left(-\bar\lambda_i^k-\log(a_k \tau_q^2)+2\right)\\
&-\frac{64\pi}{h_q^k(0)}\tau_q^{-2}e^{-\bar\lambda_i^k}
+O(e^{-(2+\epsilon_0)\bar\lambda_i^k})
\end{aligned}
\end{equation}
Now we evaluate $\int_{\Omega_l\setminus B(q,\tau_q)}\rho_kHe^{ \rm w_i^k}d\mu $ 
\begin{equation}
\label{outside-int}
\begin{aligned}
&\int_{\Omega_q\setminus B(q,\tau_q)}\rho_kHe^{ \rm w_i^k}d\mu=\int_{\Omega_q\setminus B_{\tau_q}}h_q^ke^{\bar w_i^k+\Phi}dx+O(e^{-(2+\epsilon_0)\bar \lambda_i^k})\\
&= \int_{\Omega_q\setminus B_{\tau_q}}h_q^ke^{-\bar\lambda_i^k}\frac{64}{(h_q^k(0))^2}e^{-G_q^*(q)+\Phi}dx+B_k(\tau,h)e^{-2\bar  \lambda_i^k}+O(e^{-(2+\epsilon_0)\bar \lambda_i^k})\\
&=\frac{64e^{-\bar\lambda_i^k}}{h_q^k(0)}\int_{\Omega_q\setminus B_{\tau_q}}e^{\Phi_q(x)}dx +B_k(\tau,h)e^{-2\bar  \lambda_i^k}+O(e^{-(2+\epsilon_0)\bar \lambda_i^k})
\end{aligned}
\end{equation} 
where 
\[\Phi_q(x)=\Phi(x)-G_q^*(q)+\log \frac{h_q^k(x)}{h_q^k(0)}.\]
Now we evaluate 
\begin{align*}
&\int_{\Omega_q\setminus B(q,\tau_q)}e^{\Phi_q(x)}d\mu\\
=&\int_{\Omega_q\setminus B(q,\delta)}e^{\Phi_q(x)}d\mu+\int_{B_{\delta}\setminus B_{\tau_q}}e^{\Phi_q(x)}d\mu\\
=&\int_{\Omega_q\setminus B(q,\delta)}e^{\Phi_q(x)}d\mu+
\int_{B_{\delta}\setminus B_{\tau_q}}\frac{h_q^k}{h_q^k(0)}|x|^{-4}e^{G_q^*(x)-G_q^*(q)}dx\\
=&\int_{\Omega_q\setminus B(q,\delta)}e^{\Phi_q(x)}d\mu+\pi(\tau_q^{-2}-\delta^{-2}) +\frac{\pi}2\Delta \log h_q^k(0) \log \frac{\delta}{\tau_q}+o(\tau^{\epsilon_0}).
\end{align*}
where $\epsilon_0>0$ is a small positive number. Then we see that 
\begin{equation}
\label{local-energy}
\begin{aligned}
&\int_{\Omega_q}\rho_kHe^{\rm w_i^k}d\mu\\
=&8\pi-\frac{16\pi}{h_q^k(0)}e^{-\bar\lambda_i^k}\Delta \log h_q^k(0) \left(-\bar\lambda_i^k-\log (a_k\tau_q^2)+2\right) \\
&+\frac{64e^{-\bar\lambda_i^k}}{h_q^k(0)} \left (\int_{\Omega_q\setminus B(q,\tau_q)}e^{\Phi_q(x)}d\mu-\pi \tau_q^{-2}+o(\tau^{\epsilon_0})\right)+B_k(\tau,h)e^{-2\bar  \lambda_i^k}+O(e^{-(2+\epsilon_0)\bar \lambda_i^k}). 
\end{aligned}
\end{equation}
Putting the estimates on regular points and singular sources together, we have
\begin{align*}
&\rho_k-\rho_*\\
=&\frac{16\pi}{h_q^k(0)^2e^{G_q^*(q)}}e^{-\bar\lambda_i^k}L(\mathbf{p})
\cdot \left(\bar\lambda_i^k+\log (a_k\tau_q^2)-2\right)\\
&+\frac{64 e^{-\bar\lambda_i^k}}{h_q^k(0)}\left(\int_{\Omega_q\setminus B(q,\tau_q)}e^{\Phi_q(x,q)}d\mu-\frac{\pi}{\tau_q^2}+o(\tau^{\epsilon_0})\right )  \\
&+64e^{-\tilde\lambda_i^k}e^{-G_p^*(p)}\frac{(1+\beta)^4}{h_p^k(0)^2}\left(\int_{\Omega_p\setminus B(p,\tau)}He^{\Phi}d\mu-\frac{\pi h_p^k(0)e^{G_p^*(p)}}{(1+\beta)\tau^{2+2\beta}}\right )\\
&+o(\tau^{\epsilon_0})e^{-\tilde\lambda_i^k}+B_k(\tau,h)e^{-2\bar  \lambda_i^k}+O(e^{-(2+\epsilon_0)\bar\lambda_i^k}), 
\end{align*}
Using the definition of $D(\mathbf{p})$ and relation between $\bar\lambda_i^k$ and $\tilde \lambda_i^k$, we rewrite the above as
\begin{equation}\label{combine-2}
\begin{aligned}
&\rho_k-\rho_*\\
=&\frac{16\pi}{h_q^k(0)^2e^{G_q^*(q)}}e^{-\bar\lambda_i^k}L(\mathbf{p})
\cdot \left(\bar\lambda_i^k+\log ( h_q^k(0)^2e^{G_q^*(q)}\tau^2)-2\right)  \\
&+\frac{64 e^{-\tilde\lambda_i^k}}{h_p^k(0)^2}D(\mathbf{p})+o(\tau^{\epsilon_0})e^{-\tilde\lambda_i^k}+B_k(\tau,h)e^{-2\bar  \lambda_i^k}+O(e^{-(2+\epsilon_0)\bar\lambda_i^k}), 
\end{aligned}
\end{equation}
Theorem \ref{combine-rho} is established. 

\qed
 
\medskip

We now prove that $b_0=0$. By contradiction, assume that $b_0\neq 0$, then we claim that first,
\begin{lem}\label{small-sigma}
$\sigma_k=O(\varepsilon_k^{2+\epsilon_0})$ for some $\epsilon_0>0$.
\end{lem}
\noindent{\bf Proof of Lemma \ref{small-sigma}:} In view of the definition of $\sigma_k$, \eqref{linim-2} and \eqref{linim-3}, we know that, as far as $b_0\neq 0$,
\begin{equation}\label{sigma:eq}
\sigma_k\sim |\bar \lambda_1^k-\bar \lambda_2^k|,
\end{equation}
for any $k$ large enough. The same inequality also holds for $\tilde \lambda_i^k$.
 The proof  relies on the assumption that
$\rho_1^k=\rho_2^k=\rho_k$, by looking at the estimete for $\rho_k$ in Theorem \ref{combine-rho}, we have, for some constants $C_1,C_2\neq 0$,
\begin{align*}
    \rho_i^k-\rho_{*}=&C_1\Delta \log \bar h_k(q_i^k)(\bar\lambda_i^k+\log (h_q^k(0)^2e^{G_q^*(q)}\tau^2)-2)e^{-\bar \lambda_i^k}\\
    &+C_2e^{-\bar \lambda_i^k}(D(\mathbf{p})+o(\delta))+B_k(\tau,h)e^{-2\bar  \lambda_i^k}+O(e^{-(2+\epsilon_0)\bar \lambda_i^k}).
\end{align*}
Recall that $|q_1^k-q_2^k|=o(e^{-\bar\lambda_1^k})$, which together with \eqref{same-ord} implies 
\begin{align*}
    0=&C_1 e^{-\bar \lambda_1^k} \Delta \log \bar h_q^k(q_1^k)\left( (\bar\lambda_1^k-\bar\lambda_2^k) +o(1)(\bar\lambda_1^k-\bar\lambda_2^k)\right)\\
    &+C_2e^{-\bar \lambda_1^k}(\bar\lambda_1^k-\bar\lambda_2^k)(D(\mathbf{p})+o(\delta))+O(e^{-(2+\epsilon_0)\bar \lambda_1^k}).
\end{align*}
Consequently, either $L(\mathbf{p})\neq 0$ (equivalent to $\Delta \log \tilde h^k_q(q)\neq 0$) or $L(\mathbf{p})=0$ but $D(\mathbf{p})\neq 0$, we always have 
\[\sigma_k\sim |\bar \lambda_1^k-\bar \lambda_2^k|\le C e^{-(1+\epsilon_0)\bar \lambda_1^k}.\]
Lemma \ref{small-sigma} is proved. 

\qed

\medskip

The proof of $b_0=0$ is as follows. 

\noindent{\bf Proof of $b_0=0$:} Using Lemma \ref{small-sigma}, the equation for $\xi_k$ is particularly clean now, closely resembling the linearized equation along a regular bubble. The smallness of $\sigma_k$ plays a key role in enabling us to derive a refined expansion for $\xi_k$. 
Recall that
\begin{equation*}
    \xi_k(x)=\dfrac{w_1^k(x)-w_2^k(x)}{\parallel w_1^k-w_2^k\parallel_{L^{\infty}(M)}},
\end{equation*}
and observe that, in view of Lemma \ref{small-sigma}, $c_k(x)$ can be written as follows,
\begin{equation}\label{c}
    \begin{aligned}
	c_k(x)&=\int_0^1e^{tw_1^k+(1-t)w_2^k}dt=e^{w_2^k}\left(1+\frac 12(w_1^k-w_2^k)+O(\sigma_k^2)\right),\\
	&=e^{w_2^k}(1+O(\varepsilon_k^{2+\epsilon_0})), 
    \end{aligned}
\end{equation}
which in view of \eqref{ceqtc} obviously holds for $\tilde c_k$ as well.
Therefore, in local coordinates around $p$, with $0=x(p)$, $x\in B_\tau=B_\tau(0)$, $\xi_k$ satisfies the following linearized equation,
\begin{equation*}
    \Delta\xi_k(x)+|x|^{2\beta} h_p^k(x) \tilde c_k(x)\xi_k(x)=0\quad {{\rm in}}\quad B_{\tau},
\end{equation*}
In particular, recall that the assumption $L(\mathbf{q})\neq 0$ is equivalent to $\Delta \log h_q^k(0)\neq 0$ for $k$ large enough and that we set $\bar \lambda_{i}^k=w_i^k(q) \quad i=1,2$. First of all, observe that, since locally around $p$ the rescaled function $\xi_k(\varepsilon_ky)$ (recall $0=x(p)$) converges on compact subsets of $\mathbb{R}^2$ to a solution of the following linearized equation (\cite{bart-4,bart-4-2,wu-zhang-ccm}),
$$\Delta \xi+|y|^{2\beta} h_p(0)e^{U}\xi=0,$$
where (see \eqref{st-bub}) $U$ is the limit of the standard bubbles after scaling. 
\begin{equation}\label{lim-1.0}
    \xi(y)=b_0\frac{1-c|y|^{2\beta+2}}{1+c|y|^{2\beta+2}},\quad c=\frac{ h_p(0)}{8(1+\beta)^2}.
\end{equation}
In particular, the local convergence of $\xi_k(\varepsilon_ky)$ to $\xi(y)$ readily implies that,
\begin{equation}\label{radial-xi}
    \frac{\tilde \lambda_{1}^k-\tilde \lambda_{2}^k}{\|w_1^k-w_2^k\|_{\infty}}\to b_0, \quad \frac{\bar \lambda_1^k-\bar \lambda_2^k}{\|\rm w_1^k-\rm w_2^k\|_{\infty}}\to b_0.
\end{equation}


At this point we can use (\ref{radial-xi}) and the expansions of $u_1^k$ and $u_2^k$ to identify the leading order of the radial part of $\xi_k(\varepsilon_ky)$. Let us recall that by definition $v_i^k$ is the scaling of $u_i^k$:
\begin{equation*}
    \tilde v_i^k(y)=\tilde u_{i}^k(\varepsilon_{k}y)-\phi_{i}^k(\varepsilon_{k}y)-\tilde \lambda_{i}^k,\quad i=1,2,
\end{equation*}
where $\varepsilon_{k}=e^{-\frac{\tilde \lambda_{1}^k}{2(1+\beta)}}$. Thus, by using $U_1^k$ and $U_2^k$ to denote the leading terms of $v_1^k$ and $v_2^k$, respectively, we have,
$$U_1^k-U_2^k=\tilde \lambda_1^k-\tilde \lambda_2^k+2\log \frac{1+a_ke^{\tilde \lambda_2^k-\tilde \lambda_1^k}|y|^{2\beta+2}}{1+a_k|y|^{2\beta+2}},\quad a_k= h_p^k(0)/8(1+\beta)^2.$$
Therefore,  we have,
\begin{align*}
    U_1^k(y)-U_2^k(y)
    =~&\tilde \lambda_{1}^k-\tilde \lambda_{2}^k+2\log (1+\frac{a_k|y|^{2\beta+2}(e^{\tilde \lambda_{2}^k-\tilde \lambda_{1}^k}-1)}{1+a_k|y|^{2\beta+2}})  \\
    =~&\tilde \lambda_{1}^k-\tilde \lambda_{2}^k+2\log (1+\frac{a_k|y|^{2\beta+2}(\tilde \lambda_{2}^k-\tilde \lambda_{1}^k+O(|\tilde \lambda_{2}^k-\tilde \lambda_{1}^k|^2)}{1+a_k|y|^{2\beta+2}})\\
    =~&\tilde \lambda_{1}^k-\tilde \lambda_{2}^k+2\frac{a_k|y|^{2\beta+2}(\tilde \lambda_{2}^k-\tilde \lambda_{1}^k)}{1+a_k|y|^{2\beta+2}}
    +\frac{O(\tilde \lambda_{1}^k-\tilde \lambda_{2}^k)^2}{(1+a_k|y|^{2\beta+2})^2}\\
    =~&(\tilde \lambda_{1}^k-\tilde \lambda_{2}^k)\frac{1-a_k|y|^{2\beta+2}}{1+a_k|y|^{2\beta+2}}+\frac{O(\tilde \lambda_{1}^k-\tilde \lambda_{2}^k)^2}{(1+a_k|y|^{2\beta+2})^2}.
\end{align*}

Next, we describe the expansion of $\xi_k$ near $q$. In this case \eqref{lim-1.0} is replaced by \eqref{linim-3}. Recall that we denote $\bar q_i^k=x(q_i^k)$, $i=1,2$ the expressions in local isothermal coordinates $0=x(q)$ of the maximizers near $q$ of $u_i^k$, $i=1,2$. From \eqref{close-q} we have $|\bar q_1^k-\bar q_2^k|=o(e^{-\bar \lambda_1^k})$. In $\Omega(q,\tau)$ we denote by $\bar \varepsilon_k=e^{-\bar \lambda_2^k/2}$ the scaling factor where $\bar \lambda_1^k$ and $\bar \lambda_2^k$ are the maximums near $q$. With an abuse of notations we denote by $\bar \psi^k_\xi$ the harmonic function that, in local isothermal coordinates $0=x(q)$, encodes the oscillation of $\xi_k$ on $B_\tau$.

We claim that 
\begin{equation}\label{e-b-0}
\begin{aligned}
0=&\int_M\rho_kH_k c_k \xi_kd\mu\\
=&C_1b_0\Delta \log \tilde h_k(0)\bigg (\tilde \lambda_i^k+\log (h_q^k(0)^2e^{G_q^*(q)}\tau^2\rho_*)-2\bigg )e^{-\tilde \lambda_i^k} \\
&+C_2 b_0(D(\mathbf{p})+o(\tau))e^{-\tilde \lambda_i^k}+o(e^{-\tilde \lambda_i^k}).
\end{aligned}
\end{equation}
where $C_1,C_2\neq 0 $. 
Therefore, if $L(\mathbf{p})\neq 0$, this is equivalent to $\Delta \log \tilde h_p(0)\neq 0$. Then in this case, we obviously have $b_0=0$. If $L(\mathbf{p})=0$ but $D(\mathbf{p})\neq 0$ the same conclusion also holds. 

The key computation of \eqref{e-b-0} is carried out as follows. First if we use the Green's representation formula for $\xi_k$ we have
\begin{equation}\label{exp-xi}
\xi_k(x)=\bar \xi_k+\int_MG(x,\eta)\rho_kHe^{c_k}\xi_kd\mu.
\end{equation}
If we just use a crude pointwise estimate to evaluate the difference between two points of $\xi_k$, we obtain easily 
\[|\xi_k(x_1)-\xi_k(x_2)|\le Ce^{-\bar \lambda_1^k/2},\]
for $x_1,x_2$ away from bubbling disks. This means if we use $\phi_{p,\xi}^k$ to eliminate the oscillation of $\xi_k$ on $\partial B(p,\tau)$, the estimate of $\phi_{p,\xi}^k$ is 
\[|\phi_{p,\xi}^k(\tilde \varepsilon_ky)|\le C\tilde \varepsilon_k^2|y|.\]
A similar estimate can be obtained for $\phi_{q,\xi}^k$ around $q$, where $\phi_{q,\xi}^k$ is interpreted similarly. With the correction of $\phi_{p,\xi}^k$ and $\phi_{q,\xi}^k$ we obtain like in \cite{byz-1,zhang2}
\[\xi_k(p+\tilde \varepsilon_ky)=b_0^k\frac{1-\tilde a_k|y|^{2\beta+2}}{1+\tilde a_k|y|^{2\beta+2}}+\phi_{p,\xi}^k(\tilde \varepsilon_ky)+O(\varepsilon_k^{2+\epsilon_0}(1+|y|)^{\epsilon_0})\]

The expansion around $q$ has a leading term
\[b_0^k\frac{1-\bar a_k|y|^2}{1+a_k|y|^2}+b_1^k\frac{y_1}{1+a_k |y|^2}+b_2^k\frac{y_2}{1+a_k|y|^2}.\]
The error term after expansion is also 
\[O(\bar\varepsilon_k^{2+\epsilon_0}(1+|y|)^{\epsilon_0}).\]
Using the expansions of $\xi_k$ around $p$ and $q$, \eqref{c}, and the standard calculation,  in the evaluation of $\int_M\rho_kHe^{c_k}\xi_k$, we first obtain the leading term comes from 
\begin{equation}\label{rough-1}
    \int_{B(p,\tau)}\rho_kHc_k\xi_kd\mu, \quad \int_{B(q,\tau)}\rho_kHc_k\xi_kd\mu,
\end{equation}
which is $C_1b_0\Delta \log \tilde h_k(p)\bigg (\tilde \lambda_i^k+\log (h_q^k(0)^2e^{G_q^*(q)}\tau^2\rho_*)-2\bigg )e^{-\tilde \lambda_i^k}$. Meanwhile, the error term is $O(e^{-\tilde\lambda_i^k})$.
Obviously, if $L(\mathbf{q})\neq 0$, the integral over $B(q,\tau)$ is the leading term and one obtains a contradiction to $\int_M\rho_kHc_k\xi_kd\mu=0$ immediately. If $L(\mathbf{q})=0$, the leading term is of the order $O(e^{-\bar \lambda_1^k})$. In this case, we need to evaluate each integral more precisely. In $B(p,\tau)$ and $B(q,\tau)$, using the expansion of $\xi_k$, we obviously only need to care about the radial term of $\xi_k$, which is also the leading term. For the integration outside the disks, we need to use (\ref{exp-xi}), but using the rough estimates (\ref{rough-1}), we have 
\[\xi_k=\bar \xi_k+O(\tilde\lambda_1^ke^{-\bar \lambda_1^k}).\]
Clearly $\xi_k=-\bar b_0^k(1+o(1))$. Using this expression of $\xi_k$ in the evaluation we have, under the assumption $L(\mathbf{q})=0$,
\[\int_M\rho_k c_k\xi_kd\mu
=C_2b_0(D(\mathbf{p})+o(\tau))e^{-\bar \lambda_1^k}+O(\bar\varepsilon_k^{2+\epsilon_0}).
\]
Thus, if $L(\mathbf{p})=0$ but $D(\mathbf{p})\neq 0$, we obtain a contradiction as well. 
$b_0=0$ is established. 

\qed

\medskip

\noindent{\bf Proof of $b_1=b_2=0$:}
The argument follows a strategy analogous to that in \cite{bart-4,byz-1,wu-zhang-ccm}, with a key simplification due to the condition $\alpha_M\le 0$. Specifically, the oscillation of $\xi_k$ away from the bubbling disks is $O(\bar\varepsilon_k)$, which means the harmonic function needed to remove the oscillation around $q$ is of the order $O(\bar\varepsilon_k^2|y|)$ after scaling.  This observation eliminates the need to establish a vanishing rate for $b_0^k$, unlike the case $\alpha_M>1$ treated in \cite{byz-1}, where such estimates are essential for handling positive singular sources. However, for $\alpha_M<1$, such an estimate is not needed. The presence of a negative singular source gives no extra trouble in this part. In the original paper of Bartolucci-et-al \cite{bart-4}, a vanishing rate of $b_0^k$ is also not needed. From here, we can have a glimpse of the strength of \cite{byz-1}. Just like in \cite{byz-1}, we only need to evaluate a simple single equation to prove $b_1=b_2=0$. Since this part is similar to the corresponding part in \cite{byz-1}, we omit it. 

Up to now, we have established that $b_0=b_1=b_2=0$. Theorem \ref{mainly-case-2} then follows by adapting the standard technique similar with \cite{bart-4,bart-4-2,byz-1,wu-zhang-ccm}.

\qed

\end{document}